\documentclass[twoside, 12pt]{article}

\usepackage{amsfonts}

\usepackage{indentfirst}

\begin{document} 

\title{\bf Questions on meromorphic functions and complex differential equations}

\author{Gary G. Gundersen}

\date{\today}

\maketitle

\begin{minipage}{120mm}%
\vskip.012in

\noindent{\bf Abstract.} {Thirty research questions on meromorphic functions and complex differential equations are listed and discussed. The main purpose of this paper is to make this collection of problems available to everyone.}\\\\



\end{minipage}

\baselineskip14pt

\vskip.4in

\renewcommand{\theequation}{\thesection.\arabic{equation}}

\setcounter{equation}{0}
\setcounter{section}{1}
\section*{1. Introduction}

This paper contains a list of research questions in complex analysis, most of which I wrote down many years ago and had not looked at in a long time. Partial progress, examples, and references are given.

I am grateful for the generous input and encouragement I received on my earlier drafts of the questions, both from emails and from discussions at the Workshop on Complex Differential Equations and Value Distribution Theory, held at the University of Eastern Finland in Joensuu on June 11-13, 2015. 

I would like to thank Walter Bergweiler, Yik-Man Chiang, Alex Eremenko, Janne Heittokangas, Ilpo Laine, Jim Langley, Norbert Steinmetz, and Kazuya Tohge for their invaluable help with the preparation of this paper.

Regarding future input on this collection of problems, I would be happy to receive communications concerning solutions, partial solutions, references, etc.   

Section 2 has notation, Section 3 has questions on meromorphic and entire functions, Sections 4 and 5 have questions on second order homogeneous linear complex differential equations, and Section 6 has questions on first order nonlinear complex differential equations.     

\setcounter{equation}{0}
\setcounter{section}{2}
\section*{2. Notation}  

In this paper a {\it meromorphic function} means a function that is meromorphic in the whole complex plane. We assume the reader is familiar with the Nevanlinna theory of meromorphic functions; see [34], [66]. 

For a meromorphic function $f$, we use the following notation.\\

Let $\rho(f)$ denote the order of $f$.\\

Let $\mu(f)$ denote the lower order of $f$.\\

Let $\lambda(f)$ denote the exponent of convergence of the sequence of zeros of $f$, counted according to the multiplicities.\\

Let $\overline{\lambda}(f)$ denote the exponent of convergence of the sequence of zeros of $f$, where each zero is counted once.

\setcounter{equation}{0}
\setcounter{section}{3}
\section*{3. Questions on meromorphic functions}

This section has questions on meromorphic and entire functions.\\

{\bf Question 3.1.} {\it Does there exist an entire function that possesses an infinite number of positive real zeros and no other zeros and an infinite number of purely imaginary one-points and no other one-points?}\\

If such an entire function exists, then the function would have growth of at most order two, mean type; see [12], [52].\\

{\bf Question 3.2.} {\it Does there exist an entire function that possesses an infinite number of positive real zeros and no other zeros, and an infinite number of one-points such that all the one-points lie on a finite number of half-lines that start from the origin, where an infinite number of the one-points are not positive real numbers?}\\

If such an entire function $f$ exists, then the zeros and one-points of $f$ would all lie on two or more half-lines that start from the origin, and $f$ would have finite order; see [12], [52].\\

The next four questions concern the {\it Fermat type equation}
\begin{equation}
f^n + g^n + h^n = 1,
\end{equation}
where $f, g, h$ are meromorphic functions and $n$ is a positive integer. For surveys on Questions 3.3-3.6, see [26], [29], [35].\\

{\bf Question 3.3.} {\it When either $n = 7$ or $n = 8$, does there exist nonconstant meromorphic solutions $f, g, h$ of (3.1)?}\\

The answer is no when $n \geq 9$; see [17], [18], [29], [35], [39], [41]. The answer is yes when $1 \leq n \leq 6$; see [18], [23], [25], [26], [29], [31], [51], [57].

Ishizaki [39] made partial progress on Question 3.3 by showing that if $f, g, h$ are nonconstant meromorphic functions that satisfy (3.1) when $n = 8$, then there must exist a small function $a(z)$ with respect to $f, g$ and $h$ such that 
$$W(f^8, g^8, h^8) = a(z)(f(z)g(z)h(z))^6,$$
where the left side is a Wronskian.\\

{\bf Question 3.4.} {\it When either $n = 6$ or $n = 7$, does there exist nonconstant rational solutions $f, g, h$ of (3.1)?}\\

The answer is no when $n \geq 8$; see [29], [35]. The answer is yes when $1 \leq n \leq 5$; see [26], [29], [31], [46], [51].\\

{\bf Question 3.5.} {\it When $n = 6$, does there exist nonconstant entire solutions $f, g, h$ of (3.1)?}\\

The answer is no when $n \geq 7$; see [17], [18], [29], [35], [39], [41], [64]. The answer is yes when $1 \leq n \leq 5$; see [26], [29], [31], [57].

Ishizaki [39] made partial progress on Question 3.5 by showing that if $f, g, h$ are nonconstant entire functions that satisfy (3.1) when $n = 6$, then there must exist a small function $b(z)$ with respect to $f, g$ and $h$ such that
$$W(f^6, g^6, h^6) = b(z)(f(z)g(z)h(z))^4.$$

{\bf Question 3.6.} {\it When either $n = 4$ or $n = 5$, does there exist nonconstant polynomial solutions $f, g, h$ of (3.1)?}\\

The answer is no when $n \geq 6$; see [29], [35], [51]. The answer is yes when $1 \leq n \leq 3$; see [26], [29], [46].\\

Next, recall that two nonconstant meromorphic functions $f$ and $g$ {\it share a value} $c$ in the extended plane provided that $f(z) = c$ if and only if $g(z) = c$. A value can be shared CM (counting multiplicities) or IM (ignoring multiplicities).\\

{\bf Question 3.7.} [20] {\it If $f$ and $g$ are nonconstant meromorphic functions that share three values IM and share a fourth value CM, then do $f$ and $g$ share all four values CM?}\\

The origin of Question 3.7 is the classical four point theorem of Nevanlinna [50] which states that two nonconstant meromorphic functions that share four values CM must be M\"obius transformations of each other and two of the shared values must be Picard values of the functions. It is known [20], [47] that if two nonconstant meromorphic functions share two values IM and share two other values CM, then the functions share all four values CM. In the other direction, the functions [19]
\begin{equation}
f(z) = \frac{e^z+1}{(e^z - 1)^2} \quad \hbox{and} \quad g(z) = \frac{(e^z + 1)^2}{8(e^z - 1)}
\end{equation}
share the four values $0, \infty, 1, -1/8$, where each of these four values is not shared CM. For other examples of two nonconstant meromorphic functions that share four values where each of the four values is not shared CM, see [54], [61]. 

Question 3.7 is the remaining ``gap'' between the known examples and the above ``2 CM + 2 IM = 4 CM'' theorem. 

Yi and Li [67] showed that the answer to Question 3.7 is yes when either of the following conditions is satisfied: (a) $f$ and $g$ have equal finite order that is not a positive integer, or (b) $f$ or $g$ does not have normal growth. For other partial progress on Question 3.7, see the surveys of Mues [49] and Steinmetz [62] and the references therein.\\

In the case of entire functions, Question 3.7 takes the following form.\\

{\bf Question 3.8.} {\it If $f$ and $g$ are nonconstant entire functions that share three finite values IM, then do $f$ and $g$ share all three finite values CM?}\\

Mues [48] made partial progress on Question 3.8 by showing that two nonconstant entire functions cannot share three finite values DM (different multiplicities). Here, $f$ and $g$ are said to share a value $c$ DM when the following property holds: Whenever $z_0$ is a c-point of both $f$ and $g$, the multiplicity of the c-point of $f$ is not equal to the multiplicity of the $c$-point of $g$.\\ 

Next, recall that two nonconstant meromorphic functions $f$ and $g$ {\it share a pair of values} $(a, b)$, where $a$ and $b$ are in the extended plane, provided that $f(z) = a$ if and only if $g(z) = b$. A pair of values can be shared CM or IM.\\

{\bf Question 3.9.} {\it If $f$ and $g$ are nonconstant meromorphic functions that share two pairs of values CM and share three other pairs of values IM, then do $f$ and $g$ share all five pairs CM?}\\

The conclusion of Question 3.9 implies that $f$ is a M\"obius transformation of $g$ because it is known [9], [65, p. 329] that if two nonconstant meromorphic functions share three pairs of values CM and share a fourth pair of values IM, then the functions share all four pairs CM and the functions are M\"obius transformations of each other. In the other direction, Reinders [55] observed that the functions $f$ and $g$ in (3.2) also share the pair of values $(-1/2, 1/4)$ CM, which means that $f$ and $g$ share the five pairs of values
$$(0, 0), (\infty, \infty), (1, 1), (-1/8, -1/8), (-1/2, 1/4),$$
where $(-1/2, 1/4)$ is shared CM and $(0, 0), (\infty, \infty), (1, 1), (-1/8, -1/8)$ are not shared CM. 

Thus, for five shared pairs, Question 3.9 is the remaining ``gap'' between this example and the known theorem. Hu, Li, and Yang  conjectured [38, p. 150] that the answer to Question 3.9 is yes.\\

The next question asks whether the (3.2) example is essentially the only nontrivial example of two nonconstant meromorphic functions that share five pairs.\\  

{\bf Question 3.10.} [63] {\it If $u$ and $w$ are nonconstant meromorphic functions that share five pairs of values where at least one of the five pairs is not shared CM, then does there exist M\"obius transformations $S$ and $T$, and a nonconstant entire function $h$, such that $u$ and $w$ have the forms $S(f(h))$ and $T(g(h))$, where $f$ and $g$ are the functions in (3.2)?}\\ 

Properties of meromorphic functions that share five pairs of values lead naturally to Question 3.10; see [28], [63]. If the answer to Question 3.10 is yes, then the answer to Question 3.9 would be yes. \\

The next four questions concern estimates for the modulus of the logarithmic derivative of a transcendental entire function on a curve that goes from a finite point to $\infty$.\\

{\bf Question 3.11.} {\it Let $f$ be a transcendental entire function. Does there exist a real constant $\alpha$ (positive or negative) and a curve $C$ that goes from a finite point to $\infty$, such that}
$$|f^{\prime}(z)/f(z)| \geq |z|^{\alpha}  \qquad \hbox{for all} \; z \in C \; \hbox{?}$$

\vspace{10pt}

{\bf Question 3.12.} [W. Bergweiler]  {\it Let $f$ be a transcendental entire function whose lower order $\mu(f)$ satisfies $\mu(f) > 1$. Does there exist a curve $C$ that goes from a finite point to $\infty$ such that}
$$\frac{f^{\prime}(z)}{f(z)} \to \infty  \qquad \hbox{as} \; z \to \infty \; \hbox{on} \; C \; \hbox{?}$$ 

If we were to also require that the curve $C$ in Question 3.12 lies within some arbitrarily chosen sector $S = \{z: \theta_1 < \arg z < \theta_2\},$ then the answer is no by the following example of Bergweiler: The entire function
$$f(z) = \frac{e^{-z^2} - 1}{z^2}$$
satisfies $\mu(f) = 2$ and
$$|f^{\prime}(z)/f(z)| = \frac{2 + o(1)}{|z|}$$
as $z \to \infty$ in an $\varepsilon$-sector about the positive real axis.\\

{\bf Question 3.13.} {\it Let $f$ be a transcendental entire function whose lower order $\mu(f)$ satisfies $\mu(f) > 1$. Does there exist a positive constant $\alpha$ and a curve $C$ that goes from a finite point to $\infty$, such that}
$$|f^{\prime}(z)/f(z)| \geq |z|^{\alpha}  \qquad \hbox{for all} \; z \in C \; \hbox{?}$$

If the answer to Question 3.13 is yes, then the answer to Question 3.12 would be yes.

Questions 3.11-3.13 concern lower estimates of the modulus of the logarithmic derivative of a transcendental entire function on a curve that goes from a finite point to $\infty$. For a lower estimate of the modulus of the logarithmic derivative of a transcendental entire function of finite order on a sizeable unbounded set in the plane, see [45] and the references therein.\\

{\bf Question 3.14.} {\it Let $f$ be a transcendental entire function such that for some real constant $\alpha$ satisfying $\alpha > 1$ and some curve $C$ that goes from a finite point to $\infty$, we have
$$|f^{\prime}(z)/f(z)| = O(|z|^{-\alpha}) \qquad \hbox{as} \; z \to \infty \; \hbox{on} \; C.$$
Does this imply that that there exists a complex constant $w_0$ such that $f(z) \to w_0$ as $z \to \infty$ on $C$?}\\

If the curve $C$ is a ray from the origin, then the answer is yes, see [22, Lemma 6].\\

{\bf Question 3.15.} [32, pp. 595-596], [59] {\it Let $f$ be an entire function of infinite order. Does there exist a curve $C$ that goes from a finite point to $\infty$ such that $f^{(k)}(z) \to \infty$ as $z \to \infty$ along $C$ for all $k \geq 0$?}\\

For transcendental entire functions of finite order, Langley [43] showed that the answer is yes.

\setcounter{equation}{0}
\setcounter{section}{4}
\section*{4. Questions on $f^{\prime\prime} + A(z)f = 0$}

Consider the second order homogeneous linear differential equation
\begin{equation}
f^{\prime\prime} + A(z)f = 0,
\end{equation}
where $A(z)$ is an entire function such that $A(z) \not\equiv 0$. It is well known that every solution $f$ of (4.1) is an entire function.

The questions in this section originated from papers in the 1980s. In this regard, we mention that the well known Bank-Laine Conjecture originated from the 1982 paper [2], and Bergweiler and Eremenko [7] recently solved this question.  

We use the notation for $\rho$, $\lambda$, and $\overline{\lambda}$ in Section 2.\\

{\bf Question 4.1.} {\it Let $f$ and $g$ be linearly independent solutions of (4.1), where $\overline{\lambda}(A) < \rho(A) < \infty$. Does it follow that $\lambda(fg) = \infty$?}\\

If we replace $\overline{\lambda}(A)$ with $\lambda(A)$ in Question 4.1, then a result of Bank, Laine, and Langley [5, Corollary 1] shows that the answer is yes. If in Question 4.1, we replace $\overline{\lambda}(A) < \rho(A) < \infty$ with $0 < \rho(A) < \infty$ and $\delta(0, A) = 1$, then Chiang [10] showed that the answer is yes. Here, $\delta(0, A)$ is the Nevanlinna deficiency of 0 with respect to $A(z)$. 

Under the hypothesis of Question 4.1, Bank and Laine showed that $\lambda(f) \geq \rho(A)$ for every nontrivial solution of $f$ of (4.1); see [3], [40, p. 80].\\ 

{\bf Question 4.2.} {\it Is it possible to characterize the transcendental entire functions $A(z)$ such that equation (4.1) possesses two linearly independent solutions $f, g$ satisfying $\lambda(fg) < \infty$?}\\ 

For an illustration of Question 4.2, the linearly independent functions [2]
$$f(z) = \exp((e^z - z)/2) \quad \hbox{and} \quad g(z) = \exp(-(e^z + z)/2)$$
satisfy the equation
\begin{equation}
f^{\prime\prime} - \frac{e^{2z} + 1}{4}f = 0.
\end{equation}
See also (4.3) below.\\

The next question is a particular case of Question 4.2.\\

{\bf Question 4.3.} {\it Is it possible to characterize the nonconstant periodic entire functions $A(z)$ of period $\omega$ and rational in $e^{\alpha z}$, where $\alpha = 2\pi i \omega^{-1}$, such that equation (4.1) possesses two linearly independent solutions $f, g$ satisfying $\lambda(fg) < \infty$?}\\ 

Bank and Laine [4] gave representations for every entire function $f$ satisfying $\lambda(f) < \infty$, such that $f$ is a solution of (4.1) where $A(z)$ satisfies the hypothesis in Question 4.3. 

Equation (4.2) is an example illustrating Question 4.3. Other examples are equations of the form
\begin{equation}
f^{\prime\prime} + (e^z - q^2/16) f = 0,
\end{equation}
where $q = 1, 3, 5, \cdots$, which possess linearly independent solutions $f$ and $g$ satisfying $\lambda(fg) = 0$ when $q = 1$ and $\lambda(fg) = 1$ otherwise; see [5], [40, p. 107]. For more examples and partial progress concerning Question 4.3, see [1], [4], [11] and the references therein.\\

{\bf Question 4.4.} [6] {\it If $\lambda(A) < \rho(A)$, then does every nontrivial solution $f$ of (4.1) satisfy $\lambda(f) = \infty$?}\\

Bank and Langley asked this question about a more general equation than (4.1), see [6, p. 457]. 

They proved that [6] if $A(z)$ is a transcendental entire function of finite order with finitely many zeros, then every nontrivial solution $f$ of 
$$f^{(k)} + A(z)f = 0,$$ 
where $k \geq 2$, satisfies $\lambda(f) = \infty$. This gives a partial result on Question 4.4 and answered a question that was posed in [5].\\
 
{\bf Question 4.5.} {\it Is it possible to characterize the transcendental entire functions $A(z)$ where (4.1) admits a nontrivial solution $f$ with $\lambda(f) < \infty$, such that for every solution $g$ of (4.1) that is linearly independent with $f$, we have $\lambda(g) = \infty$?}\\

For an example illustrating Question 4.5, $f(z) = \exp(e^z)$ satisfies
\begin{equation}
f^{\prime\prime} - (e^z + e^{2z}) f = 0,
\end{equation}
and for every solution $g$ of (4.4) that is linearly independent with $f$, we have $\lambda(g) = \infty$; see [5, Corollary 1].

Bank and Laine showed that if $f$ and $g$ are linearly independent solutions of (4.1), where $A(z)$ is transcendental, such that $\lambda(fg) < \infty$, then $\lambda(h) = \infty$ holds for every nontrivial solution $h$ of (4.1) that is not a constant multiple of either $f$ or $g$; see [3, Corollary 7], [40, p. 80].\\ 

Question 4.6 below is a slight modification of a question posed by Hellerstein and Rossi [8, Problem 2.71]. Bank [1983, unpublished] also posed [8, Problem 2.71].\\ 

{\bf Question 4.6.} {\it What is the most explicit form that can be found for those nonconstant polynomials $A(z)$ for which (4.1) admits a solution with an infinite number of real zeros and at most a finite number of non-real zeros?}\\

Eremenko and Merenkov [15] showed that if $n$ is any non-negative integer satisfying $n \not= 2, 6, 10, \cdots$, then there exists a polynomial $A(z)$ of degree $n$ such that equation (4.1) admits a solution $f$ that has an infinite number of real zeros and no other zeros. On the other hand, if $f$ is a nontrivial solution of (4.1) where $A(z)$ is a polynomial of degree $n$ such that $n \in \{2, 6, 10, \cdots \}$, then [21], [37] either $f$ has only finitely many zeros or the exponent of convergence of the non-real zero sequence of $f$ is $(n + 2)/2$. 

For further partial progress on Question 4.6, see [14], [24], [60] and the references therein.

\setcounter{equation}{0}
\setcounter{section}{5}
\section*{5. Questions on $f^{\prime\prime} + A(z)f^{\prime} + B(z)f = 0$}

Consider the second order homogeneous linear differential equation
\begin{equation}
f^{\prime\prime} + A(z)f^{\prime} + B(z)f = 0,
\end{equation}
where $A(z)$ and $B(z)$ are entire functions such that $A(z)B(z) \not\equiv 0$. It is well known that every solution $f$ of (5.1) is an entire function.\\

{\bf Question 5.1.} {\it If $\lambda(A) < \rho(A) < \infty$ and $B(z)$ is a nonconstant polynomial, then does every nontrivial solution $f$ of (5.1) satisfy $\rho(f) = \infty$?}\\

An example illustrating Question 5.1 is: If $B(z)$ is a nonconstant polynomial, then every nontrivial solution $f$ of
$$f^{\prime\prime} + e^{-z}f^{\prime} + B(z)f= 0$$
satisfies $\rho(f) = \infty$; see [42].  

If $B(z)$ is a constant, then the answer to Question 5.1 is no, since $f(z) = e^{-z} - 1$ satisfies $f^{\prime\prime} + e^zf^{\prime} - f = 0.$ See also [4], [16]. 

If $B(z)$ is transcendental, then the answer is also no, since $f(z) = e^{-z}$ satisfies $f^{\prime\prime} + e^zf^{\prime} + (e^z - 1)f = 0$. \\

The next question has a similar nature to Question 5.1, where $\delta(0, A)$ is the usual Nevanlinna deficiency.\\

{\bf Question 5.2.}  [J. Heittokangas] {\it If $\rho(A) < \infty$ and $\delta(0, A) > 0$, and if $B(z)$ is a nonconstant polynomial, then does every nontrivial solution $f$ of (5.1) satisfy $\rho(f) = \infty$?}\\

{\bf Question 5.3.}  {\it Does there exist an equation of the form (5.1) where $\rho(B) < \rho(A) < \infty$, that possesses two linearly independent solutions $f$ and $g$ satisfying $\lambda(fg) < \infty$?}\\

If we replace $\rho(B) < \rho(A) < \infty$ with $\rho(B) = \rho(A) < \infty$ where $A(z)$ is transcendental, then the answer is yes, since the linearly independent functions
$$f(z) = \sin z \, \exp((z + e^z)/2) \quad \hbox{and} \quad g(z) = \cos z \, \exp((z + e^z)/2)$$
satisfy the equation
$$f^{\prime\prime} - (e^z + 1) f^{\prime} + \frac{e^{2z} + 5}{4} f = 0.$$

{\bf Question 5.4.}  {\it If $\rho(B) < \rho(A) < \infty$ where $\rho(A)$ is not an integer, then does every nontrivial solution $f$ of (5.1) satisfy $\rho(f) = \infty$?}\\

If $\rho(A)$ is allowed to be a positive integer $k$ in Question 5.4, then Examples 1 and 2 in [22] show that the answer is no. For instance, by taking $Q(z) = z$ in [22, Example 1] gives the following example for $k = 1$: If $B(z)$ is any entire function satisfying $\rho(B) < 1$ and $B \not\equiv 0$, then $f(z) = e^z + 1$ satisfies the equation
$$f^{\prime\prime} - (1 + B(z) + B(z)e^{-z})f^{\prime} + B(z)f = 0,$$
where $\rho(B) < \rho(A) = \rho(f) = 1$.\\

Next, recall [36], [58, p. 6] that an entire function $f$ is said to be of {\it completely regular growth} provided that $f$ has order $\rho \in (0, \infty)$ and is of finite type, and if for each $\theta \in [0, 2\pi)$,
$$\frac{\log |f(re^{i\theta})|}{r^{\rho}} \to h_f(\theta)$$
as $r \to \infty$ outside a possible exceptional set $E$ of zero upper density where $E$ is the same set for each $\theta$. Here, $h_f(\theta)$ is the Phragm\'en-Lindel\"of indicator function of $f$.

The following question is a particular case of a general problem of Gol'dberg, Ostrovski\u{i}, and Petrenko.\\

{\bf Question 5.5.} [33, p. 300], [36] {\it Let $f$ be a transcendental solution of (5.1), where $A(z)$ and $B(z)$ have completely regular growth, such that $\rho(f) < \infty$. Does $f$ have completely regular growth?}\\

Gol'dberg [33, p. 300] gave an example which shows that the answer is no when $A(z)$ and $B(z)$ do not have completely regular growth. 

The word ``transcendental'' cannot be deleted from Question 5.5, since [36] $f(z) = z$ satisfies $f^{\prime\prime} + ze^zf^{\prime} - e^zf = 0.$

If $A(z)$ and $B(z)$ are polynomials, then [33, p. 300], [53, p. 110] a transcendental solution $f$ of (5.1) is of completely regular growth.

\setcounter{equation}{0}
\setcounter{section}{6}
\section*{6. Questions on nonlinear differential equations}

This section has questions on particular nonlinear first order complex differential equations.\\  

{\bf Question 6.1.}  {\it Consider a differential equation of the form
\begin{equation}
f^{\prime} = \frac{A_0(z) + A_1(z)f + \cdots + A_n(z)f^n}{B_0(z) + B_1(z)f + \cdots + B_m(z)f^m},
\end{equation}
where each $A_k(z), B_k(z)$ is a meromorphic function, such that (6.1) does not reduce to a Riccati or linear equation. Does there exist an equation of the form (6.1) that possesses an infinite number of distinct meromorphic solutions?}\\

It is possible for a Riccati or linear differential equation to possess an infinite number of distinct meromorphic solutions, since for every constant $c$, (i) $f(z) = \tan (z + c)$ satisfies $f^{\prime} = 1 + f^2$, (ii) $g(z) = ce^z$ satisfies $g^{\prime} = g$, and (iii) $h(z) = z^2 + c$ satisfies $h^{\prime} = 2z$.  

If all the coefficients in (6.1) are rational functions, then the answer to Question 6.1 is no, by first applying the classical Malmquist theorem [44] to obtain that every meromorphic solution of such an equation must be a rational function and then applying a general theorem of Eremenko [13] on rational solutions of first order differential equations with rational coefficients.\\ 

{\bf Question 6.2.} [30] {\it Consider a differential equation of the form 
\begin{equation}
f^{\prime} = R(e^z, f)
\end{equation}
where $R(e^z, f)$ is a rational function in $e^z$ and $f$, such that (6.2) does not reduce to a Riccati or linear equation. If $f$ is a meromorphic solution of (6.2), then does $f$ have the form  
$$f(z) = S(e^{cz}),$$
where $S(z)$ is a rational function in $z$ and $c$ is a rational number?}\\

An example [30] illustrating Question 6.2 is: 
$$f(z) = e^{-z/2} \quad \hbox{satisfies} \quad f^{\prime} = -\frac{1}{2}e^zf^3, \quad \hbox{where} \quad c = -\frac{1}{2}.$$
If we allow (6.2) to be a Riccati equation, then the answer to Question 6.2 is no, since [30]
$$f(z) = \frac{e^z}{z(e^z - 1)}$$
satisfies
$$f^{\prime} = \frac{1}{1-e^z}f + \frac{1-e^z}{e^z}f^2.$$ 
For partial positive results on Question 6.2, see [30, Theorem 2] and [56, Theorem 4.8]. For more information related to Question 6.2, see [30], [56].\\ 

The next two questions concern the first order differential equation
\begin{equation}
f^{\prime} = P_0(z) + P_1(z)f + P_2(z)f^2 + \cdots + P_n(z)f^n, \quad n \geq 3,
\end{equation}
where each $P_k(z)$ is a polynomial with $P_n(z) \not \equiv 0$. From Malmquist's theorem [44], every meromorphic solution $f$ of (6.3) must be a rational function. It is known [13], [30] that there can exist at most a finite number of distinct meromorphic (rational) solutions of (6.3). 

Let $d$ denote the number of distinct zeros of $P_n(z)$ in (6.3).\\

{\bf Question 6.3.} [27] {\it It is known that (6.3) can possess at most $M$ distinct meromorphic solutions, where
\begin{equation}
M \leq 1 + (n-1)(d^2-d+1).
\end{equation} 
Can (6.4) be improved when $d \geq 2$?}\\

The inequality (6.4) is sharp when $d = 0$ or $d = 1$, see [27, p. 5].\\

{\bf Question 6.4.} [27] {\it It is known that (6.3) can possess at most $L$ linearly independent meromorphic solutions, where 
\begin{equation}
L \leq d^2-d+2 \quad \hbox{if} \;\;\;  d \geq 3.
\end{equation}
Can (6.5) be improved?}\\

It is known that $L \leq 2$ if $0 \leq d \leq 1$ and $L \leq 3$ if $d = 2$, and these inequalities are all sharp; see [27, Section 3].

\vspace{15pt}

\noindent{\it email address:} ggunders@uno.edu

\end{document}